# ASYMPTOTIC BEHAVIOR OF WEIGHTED QUADRATIC AND CUBIC VARIATIONS OF FRACTIONAL BROWNIAN MOTION

By Ivan Nourdin

*Université Paris VI*

The present article is devoted to a fine study of the convergence of renormalized weighted quadratic and cubic variations of a fractional Brownian motion $B$ with Hurst index $H$. In the quadratic (resp. cubic) case, when $H < 1/4$ (resp. $H < 1/6$), we show by means of Malliavin calculus that the convergence holds in $L^2$ toward an explicit limit which only depends on $B$. This result is somewhat surprising when compared with the celebrated Breuer and Major theorem.

**1. Introduction and main result.** The study of single path behavior of stochastic processes is often based on the study of their power variations and there exists a very extensive literature on the subject. Recall that, a real $\kappa > 1$ being given, the $\kappa$-power variation of a process $X$, with respect to a subdivision $\pi_n = \{0 = t_{n,0} < t_{n,1} < \cdots < t_{n,n} = 1\}$ of $[0,1]$, is defined to be the sum

$$\sum_{k=0}^{n-1} |X_{t_{n,k+1}} - X_{t_{n,k}}|^\kappa.$$

For simplicity, consider from now on the case where $t_{n,k} = k/n$, for $n \in \mathbb{N}^*$ and $k \in \{0, \ldots, n\}$. In the present paper, we wish to point out some interesting phenomena when $X = B$ is a fractional Brownian motion and when the value of $\kappa$ is 2 or 3. In fact, we will also drop the absolute value (when $\kappa = 3$) and we will introduce some weights. More precisely, we will consider

$$(1.1) \qquad \sum_{k=0}^{n-1} h(B_{k/n})(\Delta B_{k/n})^\kappa, \qquad \kappa \in \{2,3\},$$

where the function $h:\mathbb{R} \to \mathbb{R}$ is assumed to be smooth enough and where $\Delta B_{k/n}$ denotes the increment $B_{(k+1)/n} - B_{k/n}$.

The analysis of the asymptotic behavior of quantities of type (1.1) is motivated, for instance, by the study of the exact rates of convergence of some approximation schemes of scalar stochastic differential equations driven by $B$ (see [5, 10] and [11]), besides, of course, the traditional applications of quadratic variations to parameter estimation problems.

Now, let us recall some known results concerning the $\kappa$-power variations (for $\kappa = 2, 3, 4, \ldots$), which are today more or less classical. First, assume that the Hurst index $H$ of $B$ is $1/2$, that is, $B$ is the *standard* Brownian motion. Let $\mu_\kappa$ denote the $\kappa$-moment of a standard Gaussian random variable $G \sim \mathcal{N}(0,1)$. By the scaling property of the Brownian motion and using the central limit theorem, it is immediate that, as $n \to \infty$:

$$(1.2) \qquad \frac{1}{\sqrt{n}} \sum_{k=0}^{n-1} [n^{\kappa/2}(\Delta B_{k/n})^\kappa - \mu_\kappa] \xrightarrow{\text{Law}} \mathcal{N}(0, \mu_{2\kappa} - \mu_\kappa^2).$$

When weights are introduced, an interesting phenomenon appears: instead of Gaussian random variables, we rather obtain *mixing* random variables as limit in (1.2). Indeed, when $\kappa$ is even, it is a very particular case of a more general result by Jacod [7] (see also [13]) that we have, as $n \to \infty$:

$$(1.3) \quad \frac{1}{\sqrt{n}} \sum_{k=0}^{n-1} h(B_{k/n})[n^{\kappa/2}(\Delta B_{k/n})^\kappa - \mu_\kappa] \xrightarrow{\text{Law}} \sqrt{\mu_{2\kappa} - \mu_\kappa^2} \int_0^1 h(B_s)\, dW_s.$$

Here, $W$ denotes another standard Brownian motion, independent of $B$. When $\kappa$ is odd, we have this time, as $n \to \infty$:

$$(1.4) \qquad \begin{aligned} &\frac{1}{\sqrt{n}} \sum_{k=0}^{n-1} h(B_{k/n})[n^{\kappa/2}(\Delta B_{k/n})^\kappa] \\ &\xrightarrow{\text{Law}} \int_0^1 h(B_s)(\sqrt{\mu_{2\kappa} - \mu_{\kappa+1}^2}\, dW_s + \mu_{\kappa+1}\, dB_s); \end{aligned}$$

see [13].

Second, assume that $H \neq 1/2$, that is, the case where the fractional Brownian motion $B$ has no independent increments anymore. Then (1.2) has been extended by Breuer and Major [1], Dobrushin and Major [3], Giraitis and Surgailis [4] or Taqqu [16]. Precisely, four cases are considered according to the evenness of $\kappa$ and the value of $H$:

• If $\kappa$ is even and if $H \in (0, 3/4)$, as $n \to \infty$,

$$(1.5) \qquad \frac{1}{\sqrt{n}} \sum_{k=0}^{n-1} [n^{\kappa H}(\Delta B_{k/n})^\kappa - \mu_\kappa] \xrightarrow{\text{Law}} \mathcal{N}(0, \sigma_{H,\kappa}^2).$$

- If $\kappa$ is even and if $H \in (3/4, 1)$, as $n \to \infty$,

$$n^{1-2H} \sum_{k=0}^{n-1} [n^{\kappa H}(\Delta B_{k/n})^\kappa - \mu_\kappa] \xrightarrow{\text{Law}} \text{``Rosenblatt r.v.''}$$

- If $\kappa$ is odd and if $H \in (0, 1/2]$, as $n \to \infty$,

(1.6) $$\frac{1}{\sqrt{n}} \sum_{k=0}^{n-1} n^{\kappa H}(\Delta B_{k/n})^\kappa \xrightarrow{\text{Law}} \mathcal{N}(0, \sigma_{H,\kappa}^2).$$

- If $\kappa$ is odd and if $H \in (1/2, 1)$, as $n \to \infty$,

$$n^{-H} \sum_{k=0}^{n-1} n^{\kappa H}(\Delta B_{k/n})^\kappa \xrightarrow{\text{Law}} \mathcal{N}(0, \sigma_{H,\kappa}^2).$$

Here, $\sigma_{H,\kappa} > 0$ denotes a constant depending only on $H$ and $\kappa$, which may be different from one formula to another one, and which can be computed explicitly. The term "Rosenblatt r.v." denotes a random variable whose distribution is the same as that of $Z$ at time one, for $Z$ the Rosenblatt process defined in [16].

Now, let us proceed with the results concerning the *weighted* power variations in the case where $H \neq 1/2$. In what follows, $h$ denotes a regular enough function such that $h$ together with its derivatives has subexponential growth. If $\kappa$ is even and $H \in (1/2, 3/4)$, then by Theorem 2 in León and Ludeña [9] (see also Corcuera, Nualart and Woerner [2] for related results on the asymptotic behavior of the $p$-variation of stochastic integrals with respect to $B$) we have, as $n \to \infty$:

(1.7) $$\frac{1}{\sqrt{n}} \sum_{k=0}^{n-1} h(B_{k/n})[n^{\kappa H}(\Delta B_{k/n})^\kappa - \mu_\kappa] \xrightarrow{\text{Law}} \sigma_{H,\kappa} \int_0^1 h(B_s)\, dW_s,$$

where, once again, $W$ denotes a standard Brownian motion independent of $B$. Thus, (1.7) shows for (1.1) a similar behavior to that observed in the standard Brownian case; compare with (1.3). In contradistinction, the asymptotic behavior of (1.1) can be completely different from (1.3) or (1.7) for other values of $H$. The first result in this direction has been observed by Gradinaru, Russo and Vallois [6] and continued in [5]. Namely, if $\kappa$ is odd and $H \in (0, 1/2)$, we have, as $n \to \infty$:

(1.8) $$n^{H-1} \sum_{k=0}^{n-1} h(B_{k/n})[n^{\kappa H}(\Delta B_{k/n})^\kappa] \xrightarrow{L^2} -\frac{\mu_{\kappa+1}}{2} \int_0^1 h'(B_s)\, ds.$$

Before giving the main result of this paper, let us make three comments. First, we stress that the limit obtained in (1.8) does not involve an independent standard Brownian motion anymore, as was the case for (1.3) or (1.7).

Second, notice that (1.8) agrees with (1.6) because, when $H \in (0, 1/2)$, we have $(\kappa + 1)H - 1 < \kappa H - 1/2$. Thus, (1.8) with $h \equiv 1$ is actually a corollary of (1.6). Third, observe that the same type of convergence as (1.8) with $H = 1/4$ had already been performed in [8], Theorem 4.1, when in (1.8) the fractional Brownian motion $B$ with Hurst index $1/4$ is replaced by an iterated Brownian motion $Z$. It is not very surprising, since this latter process is also centered, self-similar of index $1/4$ and has stationary increments. Finally, let us mention that Swanson announced in [15] that, in a joint work with Burdzy, they will prove that the same also holds for the solution of the stochastic heat equation driven by a space–time white noise.

Now, let us go back to our problem. In the sequel, we will make use of the following hypothesis on real function $h$:

($\mathbf{H_m}$) The function $h$ belongs to $\mathscr{C}^m$ and, for any $p \in (0, \infty)$ and any $0 \leq i \leq m$, we have $\sup_{t \in [0,1]} E\{|h^{(i)}(B_t)|^p\} < \infty$.

The aim of the present work is to prove the following result:

THEOREM 1.1. *Let $B$ be a fractional Brownian motion with Hurst index $H$. Then:*

1. *If $h: \mathbb{R} \to \mathbb{R}$ verifies ($\mathbf{H_4}$) and if $H \in (0, 1/4)$, we have, as $n \to \infty$:*

$$(1.9) \quad n^{2H-1} \sum_{k=0}^{n-1} h(B_{k/n})[n^{2H}(\Delta B_{k/n})^2 - 1] \xrightarrow{L^2} \tfrac{1}{4} \int_0^1 h''(B_u)\, du.$$

2. *If $h: \mathbb{R} \to \mathbb{R}$ verifies ($\mathbf{H_6}$) and if $H \in (0, 1/6)$, we have, as $n \to \infty$:*

$$(1.10) \quad \begin{aligned} & n^{3H-1} \sum_{k=0}^{n-1} [h(B_{k/n})n^{3H}(\Delta B_{k/n})^3 + \tfrac{3}{2}h'(B_{k/n})n^{-H}] \\ & \xrightarrow{L^2} -\tfrac{1}{8} \int_0^1 h'''(B_u)\, du. \end{aligned}$$

Before giving the proof of Theorem 1.1, let us roughly explain why (1.9) is only available when $H < 1/4$ [of course, the same type of argument could also be applied to understand why (1.10) is only available when $H < 1/6$]. For this purpose, let us first consider the case where $B$ is the standard Brownian motion (i.e., when $H = 1/2$). By using the independence of increments, we easily compute

$$E\left\{\sum_{k=0}^{n-1} h(B_{k/n})[n^{2H}(\Delta B_{k/n})^2 - 1]\right\} = 0$$

and

$$E\left\{\left(\sum_{k=0}^{n-1} h(B_{k/n})[n^{2H}(\Delta B_{k/n})^2 - 1]\right)^2\right\} = 2E\left\{\sum_{k=0}^{n-1} h^2(B_{k/n})\right\}$$
$$\approx 2nE\left\{\int_0^1 h^2(B_u)\,du\right\}.$$

Although these two facts are of course not sufficient to guarantee that (1.3) holds when $\kappa = 2$, they however roughly explain why it is true. Now, let us go back to the general case, that is, the case where $B$ is a fractional Brownian motion of index $H \in (0, 1/2)$. In the sequel, we will show (see Lemmas 2.2 and 2.3 for precise statements) that

$$E\left\{\sum_{k=0}^{n-1} h(B_{k/n})[n^{2H}(\Delta B_{k/n})^2 - 1]\right\} \approx \tfrac{1}{4} n^{-2H} \sum_{k=0}^{n-1} E[h''(B_{k/n})],$$

and, when $H \in (0, 1/4)$:

$$E\left\{\left(\sum_{k=0}^{n-1} h(B_{k/n})[n^{2H}(\Delta B_{k/n})^2 - 1]\right)^2\right\}$$
$$\approx \tfrac{1}{16} n^{-4H} \sum_{k,\ell=0}^{n-1} E[h''(B_{k/n})h''(B_{\ell/n})]$$
$$\approx \tfrac{1}{16} n^{2-4H} E\left\{\int\int_{[0,1]^2} h''(B_u)h''(B_v)\,du\,dv\right\}.$$

This explains the convergence (1.9). At the opposite, when $H \in (1/4, 1/2)$, one can prove that

$$E\left\{\left(\sum_{k=0}^{n-1} h(B_{k/n})[n^{2H}(\Delta B_{k/n})^2 - 1]\right)^2\right\} \approx \sigma_H^2 nE\left\{\int_0^1 h^2(B_u)\,du\right\}.$$

Thus, when $H \in (1/4, 1/2)$, the quantity $\sum_{k=0}^{n-1} h(B_{k/n})[n^{2H}(\Delta B_{k/n})^2 - 1]$ behaves as in the standard Brownian motion case, at least for the first- and second-order moments. In particular, one can expect that the following convergence holds when $H \in (1/4, 1/2)$: as $n \to \infty$,

$$(1.11) \quad \frac{1}{\sqrt{n}} \sum_{k=0}^{n-1} h(B_{k/n})[n^{2H}(\Delta B_{k/n})^2 - 1] \xrightarrow{\text{Law}} \sigma_H \int_0^1 h(B_s)\,dW_s,$$

with $W$ a standard Brownian motion independent of $B$. In fact, in the sequel of the present paper, which is a joint work with Nualart and Tudor [12], we show that (1.11) is true and we also investigate the case where $H \geq 3/4$.

Finally, let us remark that, of course, convergence (1.9) agrees with convergence (1.5). Indeed, we have $2H - 1 < -1/2$ if and only if $H < 1/4$ (it is another fact explaining the condition $H < 1/4$ in the first point of Theorem 1.1). Thus, (1.9) with $h \equiv 1$ is actually a corollary of (1.5). Similarly, (1.10) agrees with (1.5), since we have $3H - 1 < -1/2$ if and only if $H < 1/6$ (it explains the condition $H < 1/6$ in the second point of Theorem 1.1).

Now, the rest of our article is devoted to the proof of Theorem 1.1. Instead of the *pedestrian* technique performed in [6] (as their authors called it themselves), we stress the fact that we chose here a more elegant way via Malliavin calculus. It can be viewed as another novelty of this paper.

## 2. Proof of the main result.

2.1. *Notation and preliminaries.* We begin by briefly recalling some basic facts about stochastic calculus with respect to a fractional Brownian motion. One may refer to [14] for further details. Let $B = (B_t)_{t \in [0,1]}$ be a fractional Brownian motion with Hurst parameter $H \in (0, 1/2)$ defined on a probability space $(\Omega, \mathscr{A}, P)$. We mean that $B$ is a centered Gaussian process with the covariance function $\mathrm{E}(B_s B_t) = R_H(s,t)$, where

$$(2.1) \qquad R_H(s,t) = \tfrac{1}{2}(t^{2H} + s^{2H} - |t-s|^{2H}).$$

We denote by $\mathscr{E}$ the set of step $\mathbb{R}$-valued functions on $[0,1]$. Let $\mathfrak{H}$ be the Hilbert space defined as the closure of $\mathscr{E}$ with respect to the scalar product

$$\langle \mathbf{1}_{[0,t]}, \mathbf{1}_{[0,s]} \rangle_{\mathfrak{H}} = R_H(t,s).$$

We denote by $|\cdot|_{\mathfrak{H}}$ the associate norm. The mapping $\mathbf{1}_{[0,t]} \mapsto B_t$ can be extended to an isometry between $\mathfrak{H}$ and the Gaussian space $\mathcal{H}_1(B)$ associated with $B$. We denote this isometry by $\varphi \mapsto B(\varphi)$.

Let $\mathscr{S}$ be the set of all smooth cylindrical random variables, that is, of the form

$$F = f(B(\phi_1), \ldots, B(\phi_n))$$

where $n \geq 1$, $f : \mathbb{R}^n \to \mathbb{R}$ is a smooth function with compact support and $\phi_i \in \mathfrak{H}$. The Malliavin derivative of $F$ with respect to $B$ is the element of $L^2(\Omega, \mathfrak{H})$ defined by

$$D_s F = \sum_{i=1}^{n} \frac{\partial f}{\partial x_i}(B(\phi_1), \ldots, B(\phi_n)) \phi_i(s), \qquad s \in [0,1].$$

In particular $D_s B_t = \mathbf{1}_{[0,t]}(s)$. As usual, $\mathbb{D}^{1,2}$ denotes the closure of the set of smooth random variables with respect to the norm

$$\|F\|_{1,2}^2 = \mathrm{E}[F^2] + \mathrm{E}[|D.F|_{\mathfrak{H}}^2].$$

The Malliavin derivative $D$ verifies the following chain rule: if $\varphi:\mathbb{R}^n \to \mathbb{R}$ is continuously differentiable with a bounded derivative, and if $(F_i)_{i=1,\ldots,n}$ is a sequence of elements of $\mathbb{D}^{1,2}$, then $\varphi(F_1,\ldots,F_n) \in \mathbb{D}^{1,2}$ and we have, for any $s \in [0,1]$:

$$D_s\varphi(F_1,\ldots,F_n) = \sum_{i=1}^n \frac{\partial \varphi}{\partial x_i}(F_1,\ldots,F_n) D_s F_i.$$

The divergence operator $I$ is the adjoint of the derivative operator $D$. If a random variable $u \in L^2(\Omega, \mathfrak{H})$ belongs to the domain of the divergence operator, that is, if it verifies

$$|\mathrm{E}\langle DF, u\rangle_{\mathfrak{H}}| \leq c_u \|F\|_{L^2} \qquad \text{for any } F \in \mathscr{S},$$

then $I(u)$ is defined by the duality relationship

$$\mathrm{E}(FI(u)) = \mathrm{E}\langle DF, u\rangle_{\mathfrak{H}},$$

for every $F \in \mathbb{D}^{1,2}$.

2.2. *Proof of* (1.9). In this section, we assume that $H \in (0, 1/4)$. For simplicity, we note

$$\delta_{k/n} = \mathbf{1}_{[k/n,(k+1)/n]} \quad \text{and} \quad \varepsilon_{k/n} = \mathbf{1}_{[0,k/n]}.$$

Also $C$ will denote a generic constant that can be different from line to line.

We first need three lemmas. The proof of the first one follows directly from a convexity argument:

LEMMA 2.1. *For any $x \geq 0$, we have $0 \leq (x+1)^{2H} - x^{2H} \leq 1$.*

LEMMA 2.2. *For $h, g: \mathbb{R} \to \mathbb{R}$ verifying* ($\mathbf{H_2}$), *we have*

$$\sum_{k,\ell=0}^{n-1} E\{h(B_{k/n})g(B_{\ell/n})[n^{2H}(\Delta B_{k/n})^2 - 1]\}$$

(2.2)

$$= \tfrac{1}{4} n^{-2H} \sum_{k,\ell=0}^{n-1} E\{h''(B_{k/n})g(B_{\ell/n})\} + o(n^{2-2H}).$$

PROOF. For $0 \leq \ell, k \leq n-1$, we can write

$E\{h(B_{k/n})g(B_{\ell/n})n^{2H}(\Delta B_{k/n})^2\}$

$= E\{h(B_{k/n})g(B_{\ell/n})n^{2H}\Delta B_{k/n} I(\delta_{k/n})\}$

$= E\{h'(B_{k/n})g(B_{\ell/n})n^{2H}\Delta B_{k/n}\}\langle \varepsilon_{k/n}, \delta_{k/n}\rangle_{\mathfrak{H}}$

$\quad + E\{h(B_{k/n})g'(B_{\ell/n})n^{2H}\Delta B_{k/n}\}\langle \varepsilon_{\ell/n}, \delta_{k/n}\rangle_{\mathfrak{H}} + E\{h(B_{k/n})g(B_{\ell/n})\}.$

Thus,
$$n^{-2H}E\{h(B_{k/n})g(B_{\ell/n})[n^{2H}(\Delta B_{k/n})^2 - 1]\}$$
$$= E\{h'(B_{k/n})g(B_{\ell/n})I(\delta_{k/n})\}\langle\varepsilon_{k/n},\delta_{k/n}\rangle_{\mathfrak{H}}$$
(2.3)
$$+ E\{h(B_{k/n})g'(B_{\ell/n})I(\delta_{k/n})\}\langle\varepsilon_{\ell/n},\delta_{k/n}\rangle_{\mathfrak{H}}$$
$$= E\{h''(B_{k/n})g(B_{\ell/n})\}\langle\varepsilon_{k/n},\delta_{k/n}\rangle_{\mathfrak{H}}^2$$
$$+ 2E\{h'(B_{k/n})g'(B_{\ell/n})\}\langle\varepsilon_{k/n},\delta_{k/n}\rangle_{\mathfrak{H}}\langle\varepsilon_{\ell/n},\delta_{k/n}\rangle_{\mathfrak{H}}$$
$$+ E\{h(B_{k/n})g''(B_{\ell/n})\}\langle\varepsilon_{\ell/n},\delta_{k/n}\rangle_{\mathfrak{H}}^2.$$

But
(2.4)
$$\langle\varepsilon_{k/n},\delta_{k/n}\rangle_{\mathfrak{H}} = \tfrac{1}{2}n^{-2H}((k+1)^{2H} - k^{2H} - 1),$$
$$\langle\varepsilon_{\ell/n},\delta_{k/n}\rangle_{\mathfrak{H}} = \tfrac{1}{2}n^{-2H}((k+1)^{2H} - k^{2H} - |\ell-k-1|^{2H} + |\ell-k|^{2H}).$$

In particular,
(2.5)
$$|\langle\varepsilon_{k/n},\delta_{k/n}\rangle_{\mathfrak{H}}^2 - \tfrac{1}{4}n^{-4H}|$$
$$= |\tfrac{1}{4}n^{-4H}(((k+1)^{2H} - k^{2H})^2 - 2((k+1)^{2H} - k^{2H}))|$$
$$\leq \tfrac{3}{4}n^{-4H}((k+1)^{2H} - k^{2H}) \qquad \text{by Lemma 2.1}.$$

Consequently, under $(\mathbf{H_2})$:
$$n^{2H}\sum_{k,\ell=0}^{n-1}|E\{h''(B_{k/n})g(B_{\ell/n})\}(\langle\varepsilon_{k/n},\delta_{k/n}\rangle_{\mathfrak{H}}^2 - \tfrac{1}{4}n^{-4H})|$$
$$\leq Cn^{1-2H}\sum_{k=0}^{n-1}((k+1)^{2H} - k^{2H}) = Cn.$$

Similarly, using again Lemma 2.1, we deduce
$$|\langle\varepsilon_{k/n},\delta_{k/n}\rangle_{\mathfrak{H}}\langle\varepsilon_{\ell/n},\delta_{k/n}\rangle_{\mathfrak{H}}| + |\langle\varepsilon_{\ell/n},\delta_{k/n}\rangle_{\mathfrak{H}}^2|$$
$$\leq Cn^{-4H}(|(k+1)^{2H} - k^{2H}| + ||\ell-k|^{2H} - |\ell-k-1|^{2H}|).$$

Since, obviously
(2.6)
$$\sum_{k,\ell=0}^{n-1}||\ell-k|^{2H} - |\ell-k-1|^{2H}| \leq Cn^{2H+1},$$

we obtain, again under $(\mathbf{H_2})$:
$$n^{2H}\sum_{k,\ell=0}^{n-1}(|2E\{h'(B_{k/n})g'(B_{\ell/n})\}\langle\varepsilon_{k/n},\delta_{k/n}\rangle_{\mathfrak{H}}\langle\varepsilon_{\ell/n},\delta_{k/n}\rangle_{\mathfrak{H}}|$$
$$+ |E\{h(B_{k/n})g''(B_{\ell/n})\}\langle\varepsilon_{\ell/n},\delta_{k/n}\rangle_{\mathfrak{H}}^2|) \leq Cn.$$

Finally, recalling that $H < 1/4 < 1/2$, equality (2.2) follows since $n = o(n^{2-2H})$. □

LEMMA 2.3. *For $h, g : \mathbb{R} \to \mathbb{R}$ verifying* ($\mathbf{H_4}$), *we have*

(2.7)
$$\sum_{k,\ell=0}^{n-1} E\{h(B_{k/n})g(B_{\ell/n})[n^{2H}(\Delta B_{k/n})^2 - 1][n^{2H}(\Delta B_{\ell/n})^2 - 1]\}$$
$$= \tfrac{1}{16} n^{-4H} \sum_{k,\ell=0}^{n-1} E\{h''(B_{k/n})g''(B_{\ell/n})\} + o(n^{2-4H}).$$

PROOF. For $0 \leq \ell, k \leq n-1$, we can write

$$E\{h(B_{k/n})g(B_{\ell/n})[n^{2H}(\Delta B_{k/n})^2 - 1]n^{2H}(\Delta B_{\ell/n})^2\}$$
$$= E\{h(B_{k/n})g(B_{\ell/n})[n^{2H}(\Delta B_{k/n})^2 - 1]n^{2H}\Delta B_{\ell/n} I(\delta_{\ell/n})\}$$
$$= E\{h'(B_{k/n})g(B_{\ell/n})[n^{2H}(\Delta B_{k/n})^2 - 1]n^{2H}\Delta B_{\ell/n}\} \langle \varepsilon_{k/n}, \delta_{\ell/n} \rangle_{\mathfrak{H}}$$
$$+ E\{h(B_{k/n})g'(B_{\ell/n})[n^{2H}(\Delta B_{k/n})^2 - 1]n^{2H}\Delta B_{\ell/n}\} \langle \varepsilon_{\ell/n}, \delta_{\ell/n} \rangle_{\mathfrak{H}}$$
$$+ 2E\{h(B_{k/n})g(B_{\ell/n})n^{4H}\Delta B_{k/n}\Delta B_{\ell/n}\} \langle \delta_{k/n}, \delta_{\ell/n} \rangle_{\mathfrak{H}}$$
$$+ E\{h(B_{k/n})g(B_{\ell/n})[n^{2H}(\Delta B_{k/n})^2 - 1]\}.$$

Thus,

$$E\{h(B_{k/n})g(B_{\ell/n})[n^{2H}(\Delta B_{k/n})^2 - 1][n^{2H}(\Delta B_{\ell/n})^2 - 1]\}$$
$$= E\{h'(B_{k/n})g(B_{\ell/n})[n^{2H}(\Delta B_{k/n})^2 - 1]n^{2H} I(\delta_{\ell/n})\} \langle \varepsilon_{k/n}, \delta_{\ell/n} \rangle_{\mathfrak{H}}$$
$$+ E\{h(B_{k/n})g'(B_{\ell/n})[n^{2H}(\Delta B_{k/n})^2 - 1]n^{2H} I(\delta_{\ell/n})\} \langle \varepsilon_{\ell/n}, \delta_{\ell/n} \rangle_{\mathfrak{H}}$$
$$+ 2E\{h(B_{k/n})g(B_{\ell/n})n^{4H}\Delta B_{k/n} I(\delta_{\ell/n})\} \langle \delta_{k/n}, \delta_{\ell/n} \rangle_{\mathfrak{H}}$$
$$= n^{2H} E\{h''(B_{k/n})g(B_{\ell/n})[n^{2H}(\Delta B_{k/n})^2 - 1]\} \langle \varepsilon_{k/n}, \delta_{\ell/n} \rangle_{\mathfrak{H}}^2$$
$$+ 2n^{2H} E\{h'(B_{k/n})g'(B_{\ell/n})[n^{2H}(\Delta B_{k/n})^2 - 1]\}$$
$$\quad \times \langle \varepsilon_{k/n}, \delta_{\ell/n} \rangle_{\mathfrak{H}} \langle \varepsilon_{\ell/n}, \delta_{\ell/n} \rangle_{\mathfrak{H}}$$
$$+ 4n^{4H} E\{h'(B_{k/n})g(B_{\ell/n})\Delta B_{k/n}\} \langle \varepsilon_{k/n}, \delta_{\ell/n} \rangle_{\mathfrak{H}} \langle \delta_{k/n}, \delta_{\ell/n} \rangle_{\mathfrak{H}}$$
$$+ 4n^{4H} E\{h(B_{k/n})g'(B_{\ell/n})\Delta B_{k/n}\} \langle \varepsilon_{\ell/n}, \delta_{\ell/n} \rangle_{\mathfrak{H}} \langle \delta_{k/n}, \delta_{\ell/n} \rangle_{\mathfrak{H}}$$
$$+ 2n^{4H} E\{h(B_{k/n})g(B_{\ell/n})\} \langle \delta_{k/n}, \delta_{\ell/n} \rangle_{\mathfrak{H}}^2$$
$$+ n^{2H} E\{h(B_{k/n})g''(B_{\ell/n})[n^{2H}(\Delta B_{k/n})^2 - 1]\} \langle \varepsilon_{\ell/n}, \delta_{\ell/n} \rangle_{\mathfrak{H}}^2$$

$$\triangleq \sum_{i=1}^{6} A_{k,\ell,n}^{i}.$$

We claim that, for $1 \leq i \leq 5$, we have $\sum_{k,\ell=0}^{n-1} |A_{k,\ell,n}^i| = o(n^{2-4H})$. Let us first consider the case where $i=1$. By using Lemma 2.1, we deduce that

$$n^{2H} \langle \varepsilon_{k/n}, \delta_{\ell/n} \rangle_{\mathfrak{H}} = \tfrac{1}{2}((\ell+1)^{2H} - \ell^{2H} - |\ell-k+1|^{2H} + |\ell-k|^{2H})$$

is bounded. Consequently, under $(\mathbf{H_4})$:

$$|A_{k,\ell,n}^1| \leq C|\langle \varepsilon_{k/n}, \delta_{\ell/n} \rangle_{\mathfrak{H}}|.$$

As in the proof of Lemma 2.2 [see more precisely inequality (2.6)], this yields $\sum_{k,\ell=0}^{n-1} |A_{k,\ell,n}^1| \leq Cn$. Since $H < 1/4$, we finally obtain $\sum_{k,\ell=0}^{n-1} |A_{k,\ell,n}^1| = o(n^{2-4H})$. Similarly, by using the fact that $n^{2H} \langle \varepsilon_{\ell/n}, \delta_{\ell/n} \rangle_{\mathfrak{H}}$ is bounded, we prove that $\sum_{k,\ell=0}^{n-1} |A_{k,\ell,n}^2| = o(n^{2-4H})$.

Now, let us consider the case of $A_{k,\ell,n}^i$ for $i=5$, the cases where $i=3,4$ being similar. Again by Lemma 2.1, we have that

$$n^{2H} \langle \delta_{k/n}, \delta_{\ell/n} \rangle_{\mathfrak{H}} = \tfrac{1}{2}(|k-\ell+1|^{2H} + |k-\ell-1|^{2H} - 2|k-\ell|^{2H})$$

is bounded. Consequently, under $(\mathbf{H_4})$:

$$|A_{k,\ell,n}^5| \leq C(|k-\ell+1|^{2H} + |k-\ell-1|^{2H} - 2|k-\ell|^{2H}).$$

But, since $H < 1/2$, we have

$$\sum_{k,\ell=0}^{n-1} (|k-\ell+1|^{2H} + |k-\ell-1|^{2H} - 2|k-\ell|^{2H})$$

$$= \sum_{p=-\infty}^{+\infty} (|p+1|^{2H} + |p-1|^{2H} - 2|p|^{2H})$$

$$\times [(n-1) \wedge (n-1-p) - 0 \vee (-p)]$$

$$\leq Cn.$$

This yields $\sum_{k,\ell=0}^{n-1} |A_{k,\ell,n}^5| \leq Cn = o(n^{2-4H})$, again since $H < 1/4$.

It remains to consider the term with $A_{k,\ell,n}^6$. By replacing $g$ by $g''$ in identity (2.3) and by using arguments similar to those in the proof of Lemma 2.2, we can write, under $(\mathbf{H_4})$:

$$\sum_{k,\ell=0}^{n-1} A_{k,\ell,n}^6 = n^{4H} \sum_{k,\ell=0}^{n-1} E\{h''(B_{k/n})g''(B_{\ell/n})\} \langle \varepsilon_{k/n}, \delta_{k/n} \rangle_{\mathfrak{H}}^2 \langle \varepsilon_{\ell/n}, \delta_{\ell/n} \rangle_{\mathfrak{H}}^2$$

$$+ o(n^{2-4H}).$$

But, from Lemma 2.1 and equality (2.4), we deduce

(2.8)
$$|\langle \varepsilon_{k/n}, \delta_{k/n}\rangle_{\mathfrak{H}}^2 \langle \varepsilon_{\ell/n}, \delta_{\ell/n}\rangle_{\mathfrak{H}}^2 - \tfrac{1}{16}n^{-8H}|$$
$$\leq Cn^{-8H}((k+1)^{2H} - k^{2H} + (\ell+1)^{2H} - \ell^{2H}).$$

Thus,

$$n^{4H} \sum_{k,\ell=0}^{n-1} |\langle \varepsilon_{k/n}, \delta_{k/n}\rangle_{\mathfrak{H}}^2 \langle \varepsilon_{\ell/n}, \delta_{\ell/n}\rangle_{\mathfrak{H}}^2 - \tfrac{1}{16}n^{-8H}| \leq Cn^{1-2H} = o(n^{2-4H}),$$

since $H < 1/4 < 1/2$. This yields, under $(\mathbf{H_4})$:

$$\sum_{k,\ell=0}^{n-1} A_{k,\ell,n}^6 = \tfrac{1}{16}n^{-4H} \sum_{k,\ell=0}^{n-1} E\{h''(B_{k/n})g''(B_{\ell/n})\} + o(n^{2-4H}),$$

and the proof of Lemma 2.3 is done.

We are now in position to prove (1.9). Using Lemma 2.3, we have on one hand:

$$E\left\{n^{2H-1} \sum_{k=0}^{n-1} h(B_{k/n})[n^{2H}(\Delta B_{k/n})^2 - 1]\right\}^2$$

(2.9)
$$= n^{4H-2} \sum_{k,\ell=0}^{n-1} E\{h(B_{k/n})h(B_{\ell/n})$$
$$\times [n^{2H}(\Delta B_{k/n})^2 - 1][n^{2H}(\Delta B_{\ell/n})^2 - 1]\}$$
$$= \tfrac{1}{16}n^{-2} \sum_{k,\ell=0}^{n-1} E\{h''(B_{k/n})h''(B_{\ell/n})\} + o(1).$$

Using Lemma 2.2, we have on the other hand:

(2.10)
$$E\left\{n^{2H-1} \sum_{k=0}^{n-1} h(B_{k/n})[n^{2H}(\Delta B_{k/n})^2 - 1] \times \frac{1}{4n}\sum_\ell h''(B_{\ell/n})\right\}$$
$$= \sum_{k,\ell=0}^{n-1} E\{h(B_{k/n})h''(B_{\ell/n})[n^{2H}(\Delta B_{k/n})^2 - 1]\}$$
$$= \frac{1}{16}n^{-2} \sum_{k,\ell=0}^{n-1} E\{h''(B_{k/n})h''(B_{\ell/n})\} + o(1).$$

Now, we easily deduce (1.9). Indeed, thanks to (2.9)–(2.10), we obtain, by developing the square:

$$E\left\{\left(n^{2H-1}\sum_{k=0}^{n-1} h(B_{k/n})[n^{2H}(\Delta B_{k/n})^2 - 1] - \frac{1}{4n}\sum_{k=0}^{n-1} h''(B_{k/n})\right)^2\right\} \longrightarrow 0,$$

as $n \to \infty$. Since $\frac{1}{4n}\sum_{k=0}^{n-1} h''(B_{k/n}) \xrightarrow{L^2} \frac{1}{4}\int_0^1 h''(B_u)\,du$ as $n \to \infty$, we finally proved that (1.9) holds. $\square$

2.3. *Proof of* (1.10). In this section, we assume that $H \in (0, 1/6)$. We keep the same notations as in Section 2.2.

We first need two technical lemmas.

LEMMA 2.4. *For* $h, g : \mathbb{R} \to \mathbb{R}$ *verifying* $(\mathbf{H_3})$, *we have*

$$n^{3H}\sum_{k,\ell=0}^{n-1} E\{h(B_{k/n})g(B_{\ell/n})(\Delta B_{k/n})^3\}$$

$$(2.11) \qquad = -\tfrac{3}{2}n^{-H}\sum_{k,\ell=0}^{n-1} E\{h'(B_{k/n})g(B_{\ell/n})\}$$

$$-\tfrac{1}{8}n^{-3H}\sum_{k,\ell=0}^{n-1} E\{h'''(B_{k/n})g(B_{\ell/n})\} + o(n^{2-3H}).$$

PROOF. For $0 \leq \ell, k \leq n-1$, we can write

$$n^{3H}E\{h(B_{k/n})g(B_{\ell/n})(\Delta B_{k/n})^3\}$$

$$= n^{3H}E\{h(B_{k/n})g(B_{\ell/n})(\Delta B_{k/n})^2 I(\delta_{k/n})\}$$

$$= n^{3H}E\{h'(B_{k/n})g(B_{\ell/n})(\Delta B_{k/n})^2\}\langle\varepsilon_{k/n},\delta_{k/n}\rangle_{\mathfrak{H}}$$

$$\quad + n^{3H}E\{h(B_{k/n})g'(B_{\ell/n})(\Delta B_{k/n})^2\}\langle\varepsilon_{\ell/n},\delta_{k/n}\rangle_{\mathfrak{H}}$$

$$\quad + 2n^H E\{h(B_{k/n})g(B_{\ell/n})\Delta B_{k/n}\}$$

$$= n^{3H}E\{h''(B_{k/n})g(B_{\ell/n})\Delta B_{k/n}\}\langle\varepsilon_{k/n},\delta_{k/n}\rangle_{\mathfrak{H}}^2$$

$$\quad + 2n^{3H}E\{h'(B_{k/n})g'(B_{\ell/n})\Delta B_{k/n}\}\langle\varepsilon_{k/n},\delta_{k/n}\rangle_{\mathfrak{H}}\langle\varepsilon_{\ell/n},\delta_{k/n}\rangle_{\mathfrak{H}}$$

$$\quad + 3n^H E\{h'(B_{k/n})g(B_{\ell/n})\}\langle\varepsilon_{k/n},\delta_{k/n}\rangle_{\mathfrak{H}}$$

$$(2.12) \qquad + n^{3H}E\{h(B_{k/n})g''(B_{\ell/n})\Delta B_{k/n}\}\langle\varepsilon_{\ell/n},\delta_{k/n}\rangle_{\mathfrak{H}}^2$$

$$\quad + 3n^H E\{h(B_{k/n})g'(B_{\ell/n})\}\langle\varepsilon_{\ell/n},\delta_{k/n}\rangle_{\mathfrak{H}}$$

$$= n^{3H}E\{h'''(B_{k/n})g(B_{\ell/n})\}\langle\varepsilon_{k/n},\delta_{k/n}\rangle_{\mathfrak{H}}^3$$

$$\quad + 3n^{3H}E\{h''(B_{k/n})g'(B_{\ell/n})\}\langle\varepsilon_{k/n},\delta_{k/n}\rangle_{\mathfrak{H}}^2\langle\varepsilon_{\ell/n},\delta_{k/n}\rangle_{\mathfrak{H}}$$

$$\quad + 3n^{3H}E\{h'(B_{k/n})g''(B_{\ell/n})\}\langle\varepsilon_{k/n},\delta_{k/n}\rangle_{\mathfrak{H}}\langle\varepsilon_{\ell/n},\delta_{k/n}\rangle_{\mathfrak{H}}^2$$

$$\quad + 3n^H E\{h'(B_{k/n})g(B_{\ell/n})\}\langle\varepsilon_{k/n},\delta_{k/n}\rangle_{\mathfrak{H}}$$

$$+ n^{3H} E\{h(B_{k/n})g'''(B_{\ell/n})\}\langle \varepsilon_{\ell/n}, \delta_{k/n}\rangle_{\mathfrak{H}}^3$$

$$+ 3n^H E\{h(B_{k/n})g'(B_{\ell/n})\}\langle \varepsilon_{\ell/n}, \delta_{k/n}\rangle_{\mathfrak{H}} \triangleq \sum_{i=1}^{6} B_{k,\ell,n}^i.$$

We claim that $\sum_{k,\ell=0}^{n-1} |B_{k,\ell,n}^i| = o(n^{2-3H})$ for $i = 2, 3, 5, 6$. Let us first consider the cases where $i = 2$ and $i = 6$. Using Lemma 2.1 and equality (2.4), we have

$$|\langle \varepsilon_{\ell/n}, \delta_{k/n}\rangle_{\mathfrak{H}}| \leq n^{-2H}((k+1)^{2H} - k^{2H} + ||\ell - k - 1|^{2H} - |\ell - k|^{2H}|)$$

and

$$|\langle \varepsilon_{k/n}, \delta_{k/n}\rangle_{\mathfrak{H}}|^2 |\langle \varepsilon_{\ell/n}, \delta_{k/n}\rangle_{\mathfrak{H}}|$$
$$\leq Cn^{-6H}((k+1)^{2H} - k^{2H} + ||\ell - k - 1|^{2H} - |\ell - k|^{2H}|).$$

This yields, under $(\mathbf{H_3})$:

$$\sum_{k,\ell=0}^{n-1} |B_{k,\ell,n}^2| \leq Cn^{1-H} = o(n^{2-3H}) \qquad \text{since } H < 1/6 < 1/2,$$

$$\sum_{k,\ell=0}^{n-1} |B_{k,\ell,n}^6| \leq Cn^{1+H} = o(n^{2-3H}) \qquad \text{since } H < 1/6 < 1/4.$$

Similarly, we prove that $\sum_{k,\ell=0}^{n-1} |B_{k,\ell,n}^i| = o(n^{2-3H})$ for $i = 3$ and 5.

It remains to consider the terms with $B_{k,\ell,n}^1$ and $B_{k,\ell,n}^4$. From Lemma 2.1 and equality (2.4), we deduce

(2.13) $\qquad |\langle \varepsilon_{k/n}, \delta_{k/n}\rangle_{\mathfrak{H}} + \tfrac{1}{2}n^{-2H}| \leq n^{-2H}((k+1)^{2H} - k^{2H}),$

(2.14) $\qquad |\langle \varepsilon_{k/n}, \delta_{k/n}\rangle_{\mathfrak{H}}^3 + \tfrac{1}{8}n^{-6H}| \leq Cn^{-6H}((k+1)^{2H} - k^{2H}).$

Thus, since $H < 1/6$,

$$n^H \sum_{k,\ell=0}^{n-1} |\langle \varepsilon_{k/n}, \delta_{k/n}\rangle_{\mathfrak{H}} + \tfrac{1}{2}n^{-2H}| \leq n^{1+H} = o(n^{2-3H}),$$

$$n^{3H} \sum_{k,\ell=0}^{n-1} |\langle \varepsilon_{k/n}, \delta_{k/n}\rangle_{\mathfrak{H}}^3 + \tfrac{1}{8}n^{-6H}| \leq Cn^{1-H} = o(n^{2-3H}).$$

This yields, under $(\mathbf{H_3})$:

$$\sum_{k,\ell=0}^{n-1} B_{k,\ell,n}^4 = -\tfrac{3}{2}n^{-H} \sum_{k,\ell=0}^{n-1} E\{h'(B_{k/n})g(B_{\ell/n})\} + o(n^{2-3H}),$$

$$\sum_{k,\ell=0}^{n-1} B_{k,\ell,n}^1 = -\tfrac{1}{8}n^{-3H} \sum_{k,\ell=0}^{n-1} E\{h'''(B_{k/n})g(B_{\ell/n})\} + o(n^{2-3H}),$$

and the proof of Lemma 2.4 is done. $\square$

LEMMA 2.5. *For $h, g : \mathbb{R} \to \mathbb{R}$ verifying $(\mathbf{H_6})$, we have*

$$n^{6H} \sum_{k,\ell=0}^{n-1} E\{h(B_{k/n})g(B_{\ell/n})(\Delta B_{k/n})^3 (\Delta B_{\ell/n})^3\}$$

$$= \tfrac{9}{4} n^{-2H} \sum_{k,\ell=0}^{n-1} E\{h'(B_{k/n})g'(B_{\ell/n})\}$$

(2.15)
$$+ \tfrac{3}{16} n^{-4H} \sum_{k,\ell=0}^{n-1} E\{h'(B_{k/n})g'''(B_{\ell/n})\}$$

$$+ \tfrac{3}{16} n^{-4H} \sum_{k,\ell=0}^{n-1} E\{h'''(B_{k/n})g'(B_{\ell/n})\}$$

$$+ \tfrac{1}{64} n^{-6H} \sum_{k,\ell=0}^{n-1} E\{h'''(B_{k/n})g'''(B_{\ell/n})\} + o(n^{2-6H}).$$

PROOF. For $0 \leq \ell, k \leq n-1$, we can write

$n^{6H} E\{h(B_{k/n})g(B_{\ell/n})(\Delta B_{k/n})^3 (\Delta B_{\ell/n})^3\}$

$= n^{6H} E\{h(B_{k/n})g(B_{\ell/n})(\Delta B_{k/n})^3 (\Delta B_{\ell/n})^2 I(\delta_{\ell/n})\}$

$= n^{6H} E\{h'(B_{k/n})g(B_{\ell/n})(\Delta B_{k/n})^3 (\Delta B_{\ell/n})^2\} \langle \varepsilon_{k/n}, \delta_{\ell/n} \rangle_{\mathfrak{H}}$

$\quad + n^{6H} E\{h(B_{k/n})g'(B_{\ell/n})(\Delta B_{k/n})^3 (\Delta B_{\ell/n})^2\} \langle \varepsilon_{\ell/n}, \delta_{\ell/n} \rangle_{\mathfrak{H}}$

$\quad + 3n^{6H} E\{h(B_{k/n})g(B_{\ell/n})(\Delta B_{k/n})^2 (\Delta B_{\ell/n})^2\} \langle \delta_{k/n}, \delta_{\ell/n} \rangle_{\mathfrak{H}}$

$\quad + 2n^{4H} E\{h(B_{k/n})g(B_{\ell/n})(\Delta B_{k/n})^3 \Delta B_{\ell/n}\}$

$= n^{6H} E\{h''(B_{k/n})g(B_{\ell/n})(\Delta B_{k/n})^3 \Delta B_{\ell/n}\} \langle \varepsilon_{k/n}, \delta_{\ell/n} \rangle_{\mathfrak{H}}^2$

$\quad + 2n^{6H} E\{h'(B_{k/n})g'(B_{\ell/n})(\Delta B_{k/n})^3 \Delta B_{\ell/n}\} \langle \varepsilon_{k/n}, \delta_{\ell/n} \rangle_{\mathfrak{H}} \langle \varepsilon_{\ell/n}, \delta_{\ell/n} \rangle_{\mathfrak{H}}$

$\quad + 6n^{6H} E\{h'(B_{k/n})g(B_{\ell/n})(\Delta B_{k/n})^2 \Delta B_{\ell/n}\} \langle \varepsilon_{k/n}, \delta_{\ell/n} \rangle_{\mathfrak{H}} \langle \delta_{k/n}, \delta_{\ell/n} \rangle_{\mathfrak{H}}$

$\quad + 3n^{4H} E\{h'(B_{k/n})g(B_{\ell/n})(\Delta B_{k/n})^3\} \langle \varepsilon_{k/n}, \delta_{\ell/n} \rangle_{\mathfrak{H}}$

$\quad + n^{6H} E\{h(B_{k/n})g''(B_{\ell/n})(\Delta B_{k/n})^3 \Delta B_{\ell/n}\} \langle \varepsilon_{\ell/n}, \delta_{\ell/n} \rangle_{\mathfrak{H}}^2$

$\quad + 6n^{6H} E\{h(B_{k/n})g'(B_{\ell/n})(\Delta B_{k/n})^2 \Delta B_{\ell/n}\} \langle \varepsilon_{\ell/n}, \delta_{\ell/n} \rangle_{\mathfrak{H}} \langle \delta_{k/n}, \delta_{\ell/n} \rangle_{\mathfrak{H}}$

$\quad + 3n^{4H} E\{h(B_{k/n})g'(B_{\ell/n})(\Delta B_{k/n})^3\} \langle \varepsilon_{\ell/n}, \delta_{\ell/n} \rangle_{\mathfrak{H}}$

$\quad + 6n^{6H} E\{h(B_{k/n})g(B_{\ell/n}) \Delta B_{k/n} \Delta B_{\ell/n}\} \langle \delta_{k/n}, \delta_{\ell/n} \rangle_{\mathfrak{H}}^2$

$\quad + 9n^{4H} E\{h(B_{k/n})g(B_{\ell/n})(\Delta B_{k/n})^2\} \langle \delta_{k/n}, \delta_{\ell/n} \rangle_{\mathfrak{H}}$

$$= n^{6H} E\{h'''(B_{k/n})g(B_{\ell/n})(\Delta B_{k/n})^3\}\langle \varepsilon_{k/n}, \delta_{\ell/n}\rangle_{\mathfrak{H}}^3$$
$$+ 3n^{6H} E\{h''(B_{k/n})g'(B_{\ell/n})(\Delta B_{k/n})^3\}\langle \varepsilon_{k/n}, \delta_{\ell/n}\rangle_{\mathfrak{H}}^2 \langle \varepsilon_{\ell/n}, \delta_{\ell/n}\rangle_{\mathfrak{H}}$$
$$+ 9n^{6H} E\{h''(B_{k/n})g(B_{\ell/n})(\Delta B_{k/n})^2\}\langle \varepsilon_{k/n}, \delta_{\ell/n}\rangle_{\mathfrak{H}}^2 \langle \delta_{k/n}, \delta_{\ell/n}\rangle_{\mathfrak{H}}$$
$$+ 3n^{6H} E\{h'(B_{k/n})g''(B_{\ell/n})(\Delta B_{k/n})^3\}\langle \varepsilon_{k/n}, \delta_{\ell/n}\rangle_{\mathfrak{H}} \langle \varepsilon_{\ell/n}, \delta_{\ell/n}\rangle_{\mathfrak{H}}^2$$
$$+ 12n^{6H} E\{h'(B_{k/n})g'(B_{\ell/n})(\Delta B_{k/n})^2\}$$
$$\times \langle \varepsilon_{k/n}, \delta_{\ell/n}\rangle_{\mathfrak{H}} \langle \varepsilon_{\ell/n}, \delta_{\ell/n}\rangle_{\mathfrak{H}} \langle \delta_{k/n}, \delta_{\ell/n}\rangle_{\mathfrak{H}}$$
$$+ 12n^{6H} E\{h'(B_{k/n})g'(B_{\ell/n})\Delta B_{k/n}\}\langle \varepsilon_{k/n}, \delta_{\ell/n}\rangle_{\mathfrak{H}} \langle \delta_{k/n}, \delta_{\ell/n}\rangle_{\mathfrak{H}}^2$$
$$+ 3n^{4H} E\{h'(B_{k/n})g(B_{\ell/n})(\Delta B_{k/n})^3\}\langle \varepsilon_{k/n}, \delta_{\ell/n}\rangle_{\mathfrak{H}}$$
$$+ n^{6H} E\{h(B_{k/n})g'''(B_{\ell/n})(\Delta B_{k/n})^3\}\langle \varepsilon_{\ell/n}, \delta_{\ell/n}\rangle_{\mathfrak{H}}^3$$
$$+ 9n^{6H} E\{h(B_{k/n})g''(B_{\ell/n})(\Delta B_{k/n})^2\}\langle \varepsilon_{\ell/n}, \delta_{\ell/n}\rangle_{\mathfrak{H}}^2 \langle \delta_{k/n}, \delta_{\ell/n}\rangle_{\mathfrak{H}}$$
$$+ 18n^{6H} E\{h(B_{k/n})g'(B_{\ell/n})\Delta B_{k/n}\}\langle \varepsilon_{\ell/n}, \delta_{\ell/n}\rangle_{\mathfrak{H}} \langle \delta_{k/n}, \delta_{\ell/n}\rangle_{\mathfrak{H}}^2$$
$$+ 3n^{4H} E\{h(B_{k/n})g'(B_{\ell/n})(\Delta B_{k/n})^3\}\langle \varepsilon_{\ell/n}, \delta_{\ell/n}\rangle_{\mathfrak{H}}$$
$$+ 6n^{6H} E\{h'(B_{k/n})g(B_{\ell/n})\Delta B_{k/n}\}\langle \varepsilon_{k/n}, \delta_{\ell/n}\rangle_{\mathfrak{H}} \langle \delta_{k/n}, \delta_{\ell/n}\rangle_{\mathfrak{H}}^2$$
$$+ 6n^{6H} E\{h(B_{k/n})g(B_{\ell/n})(\Delta B_{k/n})^2\}\langle \delta_{k/n}, \delta_{\ell/n}\rangle_{\mathfrak{H}}^3$$
$$+ 9n^{4H} E\{h(B_{k/n})g(B_{\ell/n})(\Delta B_{k/n})^2\}\langle \delta_{k/n}, \delta_{\ell/n}\rangle_{\mathfrak{H}}.$$

To obtain (2.15), we develop the right-hand side of the previous identity in the same way as for the obtention of (2.12). Then, only the terms containing $\langle \varepsilon_{k/n}, \delta_{k/n}\rangle_{\mathfrak{H}}^{\alpha} \langle \varepsilon_{\ell/n}, \delta_{\ell/n}\rangle_{\mathfrak{H}}^{\beta}$, for $\alpha, \beta \geq 1$, have a contribution in (2.15), as we can check by using (2.5), (2.8), (2.13) and (2.14). The other terms are $o(n^{2-6H})$. Details are left to the reader. □

We are now in position to prove (1.10). Using Lemmas 2.4 and 2.5, we have on one hand

$$E\left\{\left(n^{3H-1}\sum_{k=0}^{n-1}[h(B_{k/n})n^{3H}(\Delta B_{k/n})^3 + \tfrac{3}{2}h'(B_{k/n})n^{-H}]\right)^2\right\}$$

(2.16)
$$= n^{6H-2}\sum_{k,\ell=0}^{n-1} E\{[h(B_{k/n})n^{3H}(\Delta B_{k/n})^3 + \tfrac{3}{2}h'(B_{k/n})n^{-H}]$$
$$\times [h(B_{\ell/n})n^{3H}(\Delta B_{\ell/n})^3 + \tfrac{3}{2}h'(B_{\ell/n})n^{-H}]\}$$
$$= \tfrac{1}{64}n^{-2}\sum_{k,\ell=0}^{n-1} E\{h'''(B_{k/n})h'''(B_{\ell/n})\} + o(1).$$

On the other hand, we have, by Lemma 2.4:

$$E\Bigg\{n^{3H-1}\sum_{k=0}^{n-1}\bigg[h(B_{k/n})n^{3H}(\Delta B_{k/n})^3 + \frac{3}{2}h'(B_{k/n})n^{-H}\bigg]$$
$$\times \frac{-1}{8n}\sum_{\ell} h'''(B_{\ell/n})\Bigg\}$$

(2.17)
$$= -\frac{n^{3H-2}}{8}\sum_{k,\ell=0}^{n-1} E\bigg\{h(B_{k/n})h'''(B_{\ell/n})n^{3H}(\Delta B_{k/n})^3$$
$$+ \frac{3}{2}h'(B_{k/n})h'''(B_{\ell/n})n^{-H}\bigg\}$$
$$= \frac{1}{64}n^{-2}\sum_{k,\ell=0}^{n-1} E\{h'''(B_{k/n})h'''(B_{\ell/n})\} + o(1).$$

Now, we easily deduce (1.10). Indeed, thanks to (2.16)–(2.17), we obtain, by developing the square:

$$E\Bigg\{\bigg(n^{3H-1}\sum_{k=0}^{n-1}\bigg[h(B_{k/n})n^{3H}(\Delta B_{k/n})^3 + \frac{3}{2}h'(B_{k/n})n^{-H}\bigg]$$
$$+ \frac{1}{8n}\sum_{k=0}^{n-1} h'''(B_{k/n})\bigg)^2\Bigg\} \longrightarrow 0$$
$$\text{as } n \to \infty.$$

Since $-\frac{1}{8n}\sum_{k=0}^{n-1} h'''(B_{k/n}) \xrightarrow{L^2} -\frac{1}{8}\int_0^1 h'''(B_u)\,du$ as $n \to \infty$, we finally proved that (1.10) holds.

**Acknowledgment.** I want to thank the anonymous referee whose remarks and suggestions greatly improved the presentation of my paper.

Laboratoire de Probabilités et Modèles Aléatoires
Université Paris VI
Boîte courrier 188
75252 Paris Cedex 05
France
E-mail: ivan.nourdin@upmc.fr